\documentclass[12pt]{article}
\usepackage{amsfonts,amssymb,amsmath}

\font\smallit=cmti10

\newtheorem{theorem}{Theorem}

\newtheorem{lemma}[theorem]{Lemma}
\newtheorem{proposition}[theorem]{Proposition}
\newtheorem{remark}[theorem]{Remark}

\begin{document}

\begin{center}
{\bf A NOTE ON BOOLEAN LATTICES AND FAREY SEQUENCES II} \vskip 20pt {\bf Andrey O. Matveev}\\ {\smallit
Data-Center Co., RU-620034, Ekaterinburg, P.O.~Box~5, Russian~Federation}\\ {\tt aomatveev@dc.ru,
aomatveev@hotmail.com}
\end{center}

\vskip 30pt  \vskip 30pt

\centerline{\bf Abstract}

\noindent We establish monotone bijections between subsequences of the Farey sequences and the~halfsequences
of Farey subsequences associated with elements of the Boolean lattices.

\pagestyle{myheadings}

\thispagestyle{empty} \baselineskip=15pt \vskip 30pt

\section*{\normalsize 1. Introduction}

The {\em Farey sequence of order $n$}, denoted by $\mathcal{F}_n$, is the ascending sequence of irreducible
fractions $\tfrac{h}{k}\in\mathbb{Q}$ such that $\tfrac{0}{1}\leq\tfrac{h}{k}\leq\tfrac{1}{1}$ and $1\leq
k\leq n$, see, e.g., ~[3, Chapter~27], [4, \S{}3], [5, Chapter~4], [6, Chapter~III], [11, Chapter~6], [12,
Chapter~6], [13, Sequences A006842 and~A006843], [14, Chapter~5].

The Farey sequence of order $n$ contains the subsequences
\begin{align}
\label{eq:1} \mathcal{F}_n^m&:=\left(\tfrac{h}{k}\in\mathcal{F}_n:\ h\leq m\right)\ , \intertext{for integer
$m\geq 1$, and the subsequences} \label{eq:8} \mathcal{G}_n^m&:=\left(\tfrac{h}{k}\in\mathcal{F}_n:\ k-h\leq
n-m\right)\ ,
\end{align}
for $m\leq n-1$, that inherit many familiar properties of the Farey sequences. In particular, if $n>1$ and
$m\geq n-1$, then $\mathcal{F}_n^m=\mathcal{F}_n$; if $n>1$ and $m\leq 1$, then
$\mathcal{G}_n^m=\mathcal{F}_n$.

The Farey subsequence $\mathcal{F}_n^m$ was presented in~[1], and some comments were given in~[10,
Remark~7.10].

Let $\mathbb{B}(n)$ denote the Boolean lattice of rank $n>1$, whose operation of meet is denoted by $\wedge$.
If $a$ is an element of~$\mathbb{B}(n)$, of rank $\rho(a)=:m$ such that $0<m<n$, then the integers~$n$ and
$m$ determine the subsequence
\begin{equation}
\label{eq:2}
\begin{split}
\mathcal{F}\bigl(\mathbb{B}(n),m\bigr):&=\left(\tfrac{\rho(b\wedge a)}{\gcd(\rho(b\wedge
a),\rho(b))}\!\Bigm/\!\tfrac{\rho(b)}{\gcd(\rho(b\wedge a),\rho(b))}:\ b\in\mathbb{B}(n),\ \rho(b)>0\right)\\
&=\left(\tfrac{h}{k}\in\mathcal{F}_n:\ h\leq m,\ k-h\leq n-m\right)
\end{split}
\end{equation}
of $\mathcal{F}_n$, considered in~[7, 8, 9, 10].

Notice that the map
\begin{equation}
\label{eq:3} \mathcal{F}\bigl(\mathbb{B}(n),m\bigr)\to\mathcal{F}\bigl(\mathbb{B}(n),n-m\bigr)\ ,\ \ \
\tfrac{h}{k}\mapsto\tfrac{k-h}{k}\ ,\ \ \ \left[\begin{smallmatrix}\!h\!\\
\!k\!\end{smallmatrix}\right]\mapsto\left[\begin{smallmatrix}-1&1\\0&1\end{smallmatrix}\right]\cdot
\left[\begin{smallmatrix}\!h\!\\ \!k\!\end{smallmatrix}\right]\ ,
\end{equation}
is order-reversing and bijective, by analogy with the map
\begin{equation}
\label{eq:11} \mathcal{F}_n\to\mathcal{F}_n\ ,\ \ \ \tfrac{h}{k}\mapsto\tfrac{k-h}{k}\ ,\ \ \
\left[\begin{smallmatrix}\!h\!\\
\!k\!\end{smallmatrix}\right]\mapsto\left[\begin{smallmatrix}-1&1\\0&1\end{smallmatrix}\right]\cdot
\left[\begin{smallmatrix}\!h\!\\ \!k\!\end{smallmatrix}\right]\ ,
\end{equation}
which is order-reversing and bijective as well; here $\left[\begin{smallmatrix}\!h\!\\
\!k\!\end{smallmatrix}\right]\in\mathbb{Z}^2$ is a vector presentation of the~fraction $\tfrac{h}{k}$.

We call the ascending sets
\begin{align*}
\mathcal{F}^{\leq\frac{1}{2}}\bigl(\mathbb{B}(n),m\bigr):=\left(\tfrac{h}{k}
\in\mathcal{F}\bigl(\mathbb{B}(n),m\bigr):\ \tfrac{h}{k}\leq\tfrac{1}{2}\right)\ ,\intertext{and}
\mathcal{F}^{\geq\frac{1}{2}}\bigl(\mathbb{B}(n),m\bigr):=\left(\tfrac{h}{k}\in
\mathcal{F}\bigl(\mathbb{B}(n),m\bigr):\ \tfrac{h}{k}\geq\tfrac{1}{2}\right)\ .
\end{align*}
the {\em left\/} and {\em right halfsequences\/} of the sequence $\mathcal{F}\bigl(\mathbb{B}(n),m\bigr)$,
respectively.

This work is a sequel to note~[7] which concerns the sequences $\mathcal{F}\bigl(\mathbb{B}(2m),m\bigr)$.
See~[7] for more about Boolean lattices, Farey (sub)sequences and related combinatorial identities.

The key observation is Theorem~\ref{prop:1} which in particular asserts that there is a bijection between the
sequence~$\mathcal{F}_{n-m}^m$ and the left halfsequence of $\mathcal{F}\bigl(\mathbb{B}(n),m\bigr)$, on the
one hand; there is also a bijection between the sequence $\mathcal{F}_m^{n-m}$ and the right halfsequence
of~$\mathcal{F}\bigl(\mathbb{B}(n),m\bigr)$, on the other hand.

Throughout the note, $n$ and $m$ mean integer numbers; $n$ is always greater than one. The fractions of Farey
(sub)sequences are indexed starting with zero. Terms `precedes' and `succeeds' related to pairs of fractions
always mean relations of immediate consecution. When we deal with a Farey subsequence
$\mathcal{F}\bigl(\mathbb{B}(n),m\bigr)$, we implicitly suppose $0<m<n$.

\vskip 30pt

\section*{\normalsize 2. The Farey Subsequences $\mathcal{F}_n^m$ and $\mathcal{G}_n^m$}

The connection between the sequences of the form~(\ref{eq:1}) and~(\ref{eq:8}) comes from their definitions:

\begin{lemma}
\label{prop:2}
The maps
\begin{align*}
\mathcal{F}_n^m&\to\mathcal{G}_n^{n-m}\ , & \tfrac{h}{k}&\mapsto\tfrac{k-h}{k}\ , &
\left[\begin{smallmatrix}\!h\!\\
\!k\!\end{smallmatrix}\right]&\mapsto\left[\begin{smallmatrix}-1&1\\0&1\end{smallmatrix}\right]\cdot
\left[\begin{smallmatrix}\!h\!\\ \!k\!\end{smallmatrix}\right]\ , \intertext{and}
\mathcal{G}_n^m&\to\mathcal{F}_n^{n-m}\ , & \tfrac{h}{k}&\mapsto\tfrac{k-h}{k}\ , &
\left[\begin{smallmatrix}\!h\!\\
\!k\!\end{smallmatrix}\right]&\mapsto\left[\begin{smallmatrix}-1&1\\0&1\end{smallmatrix}\right]\cdot
\left[\begin{smallmatrix}\!h\!\\ \!k\!\end{smallmatrix}\right]\ ,
\end{align*}
are order-reversing and bijective, for any $m$, $0\leq m\leq n$.
\end{lemma}

Following [2, \S{}4.15], we let $\overline{\mu}(\cdot)$ denote the M\"{o}bius function on positive integers.
For an interval of positive integers $[i,l]:=\{j:\ i\leq j\leq l\}$ and for a positive integer $h$, let
$\phi(h;[i,l]):=|\{j\in[i,l]:\ \gcd(h,j)=1\}|$.

We now describe several basic properties of the sequences $\mathcal{G}_n^m$ which are dual, in view
of~Lemma~\ref{prop:2}, to those of the sequences $\mathcal{F}_n^{n-m}$, cf.~[10, Remark~7.10].

\begin{remark}
Suppose $0\leq m<n$.
\begin{itemize}
\item[\rm(i)] In $\mathcal{G}_n^m=(g_0<g_1<\cdots<g_{|\mathcal{G}_n^m|-2}<g_{|\mathcal{G}_n^m|-1})$, we have
\begin{equation*}
g_0=\tfrac{0}{1}\ ,\ \ \ g_1=\tfrac{1}{\min\{n-m+1,n\}}\ ,\ \ \ g_{|\mathcal{G}_n^m|-2}=\tfrac{n-1}{n}\ ,\ \
\ g_{|\mathcal{G}_n^m|-1}=\tfrac{1}{1}\ .
\end{equation*}

\item[\rm(ii)]
\begin{itemize}
\item[\rm(a)]
The cardinality of the sequence $\mathcal{G}_n^m$ equals
\begin{multline*}
1+\sum_{j\in[1,n]}\phi\bigl(j;[\max\{1, j+m-n\},j]\bigr)\\=1+\sum_{j\in[1,n-m+1]}\phi\bigl(j;[1,j]\bigr)+
\sum_{j\in[n-m+2,n]}\phi\bigl(j;[j+m-n,j]\bigr)\\= 1+\sum_{d\geq
1}\overline{\mu}(d)\!\cdot\!\Bigl(\left\lfloor\tfrac{n}{d}\right\rfloor
-\tfrac{1}{2}\left\lfloor\tfrac{n-m}{d}\right\rfloor\Bigr)\!\!\cdot\!\!\left\lfloor\tfrac{n-m}{d}+1\right\rfloor\
.
\end{multline*}

\item[\rm(b)] If $g_t\in\mathcal{G}_n^m-\{\tfrac{0}{1}\}$, then
\begin{multline*}
t=\sum_{j\in[1,n]}\phi\bigl(j;[\max\{1, j+m-n\},\left\lfloor j\cdot g_t\right\rfloor]\bigr)\\=
\sum_{j\in[1,n-m+1]}\phi\bigl(j;[1,\left\lfloor j\cdot g_t\right\rfloor]\bigr)+
\sum_{j\in[n-m+2,n]}\phi\bigl(j;[j+m-n,\left\lfloor j\cdot g_t\right\rfloor]\bigr)\\= 1+ \sum_{d\geq
1}\overline{\mu}(d)\!\cdot\!\!\Biggl(\!\left\lfloor\tfrac{n-m}{d}\right\rfloor\!\!\cdot\!\!
\Bigl(\left\lfloor\tfrac{n}{d}\right\rfloor-\tfrac{1}{2}\left\lfloor\tfrac{n-m}{d}+1\right\rfloor\Bigr)\\
-\sum_{j\in[1,\lfloor n/d\rfloor]} \min\Bigl\{\left\lfloor\tfrac{n-m}{d}\right\rfloor\! , \left\lfloor
j\!\cdot\!(1-g_t)\right\rfloor\Bigr\}\!\Biggr)\ .
\end{multline*}
\end{itemize}

\item[\rm(iii)] Let $\tfrac{h}{k}\in\mathcal{G}_n^m$, $\tfrac{0}{1}<\tfrac{h}{k}<\tfrac{1}{1}$.
\begin{itemize}
\item[\rm(a)]
Let $x_0$ be the integer such that $kx_0\equiv-1\pmod{h}$ and \mbox{$m-h+1$} $\leq x_0\leq m$. Define
integers $y_0$ and $t^{\ast}$ by $y_0:=\tfrac{kx_0+1}{h}$ and~$t^{\ast}$
$:=\left\lfloor\min\left\{\tfrac{n-m+x_0-y_0}{k-h},\tfrac{n-y_0}{k}\right\}\right\rfloor$. The~fraction
$\tfrac{x_0+t^{\ast}h}{y_0+t^{\ast}k}$ precedes the~fraction $\tfrac{h}{k}$ in $\mathcal{G}_n^m$.
\item[\rm(b)]
Let $x_0$ be the integer such that $kx_0\equiv 1\pmod{h}$ and \mbox{$m-h+1$} $\leq x_0\leq m$. Define
integers $y_0$ and $t^{\ast}$ by $y_0:=\tfrac{kx_0-1}{h}$ and~$t^{\ast}$
$:=\left\lfloor\min\left\{\tfrac{n-m+x_0-y_0}{k-h},\tfrac{n-y_0}{k}\right\}\right\rfloor$.
The~fraction~$\tfrac{x_0+t^{\ast}h}{y_0+t^{\ast}k}$ succeeds the~fraction $\tfrac{h}{k}$ in
$\mathcal{G}_n^m$.
\end{itemize}

\item[\rm(iv)]
\begin{itemize}
\item[\rm(a)]
If $\tfrac{h_j}{k_j}<\tfrac{h_{j+1}}{k_{j+1}}$ are two successive fractions of $\mathcal{G}_n^m$ then
\begin{equation*}
k_j h_{j+1}-h_j k_{j+1}=1\ .
\end{equation*}
\item[\rm(b)]
If $\tfrac{h_j}{k_j}<\tfrac{h_{j+1}}{k_{j+1}}<\tfrac{h_{j+2}}{k_{j+2}}$ are three successive fractions of
$\mathcal{G}_n^m$ then
\begin{equation*}
\tfrac{h_{j+1}}{k_{j+1}}=\tfrac{h_j+h_{j+2}}{\gcd(h_j+h_{j+2},k_j+k_{j+2})}\!\!\Bigm/\!\!
\tfrac{k_j+k_{j+2}}{\gcd(h_j+h_{j+2},k_j+k_{j+2})}\ .
\end{equation*}
\item[\rm(c)]
If $\tfrac{h_j}{k_j}<\tfrac{h_{j+1}}{k_{j+1}}<\tfrac{h_{j+2}}{k_{j+2}}$ are three successive fractions of
$\mathcal{G}_n^m$ then the integers $h_j$, $k_j$, $h_{j+2}$ and $k_{j+2}$ are related in the following way:
\begin{align*}
h_j&=\left\lfloor\min\left\{\tfrac{k_{j+2}+n}{k_{j+1}},\tfrac{k_{j+2}
-h_{j+2}+n-m}{k_{j+1}-h_{j+1}}\right\}\right\rfloor h_{j+1}-h_{j+2}\ ,\\
k_j&=\left\lfloor\min\left\{\tfrac{k_{j+2}+n}{k_{j+1}},\tfrac{k_{j+2}
-h_{j+2}+n-m}{k_{j+1}-h_{j+1}}\right\}\right\rfloor k_{j+1}-k_{j+2}\ ,\\
h_{j+2}&=\left\lfloor\min\left\{\tfrac{k_j+n}{k_{j+1}},\tfrac{k_j
-h_j+n-m}{k_{j+1}-h_{j+1}}\right\}\right\rfloor h_{j+1}-h_j\ ,\\
k_{j+2}&=\left\lfloor\min\left\{\tfrac{k_j+n}{k_{j+1}},\tfrac{k_j
-h_j+n-m}{k_{j+1}-h_{j+1}}\right\}\right\rfloor k_{j+1}-k_j\ .
\end{align*}
\end{itemize}

\item[\rm(v)]
If $\tfrac{1}{k}\in\mathcal{G}_n^m$, where $n>1$, for some $k>1$, then the fraction\newline
$\tfrac{\lfloor\min\{\frac{n-m-1}{k-1},\frac{n-1}{k}\}\rfloor}{k
\lfloor\min\{\frac{n-m-1}{k-1},\frac{n-1}{k}\}\rfloor+1}$ precedes $\tfrac{1}{k}$, and the fraction
$\tfrac{\lfloor\min\{\frac{n-m+1}{k-1},\frac{n+1}{k}\}\rfloor}{k
\lfloor\min\{\frac{n-m+1}{k-1},\frac{n+1}{k}\}\rfloor-1}$ succeeds $\tfrac{1}{k}$ in $\mathcal{G}_n^m$.
\end{itemize}
\end{remark}

\vskip 30pt

\section*{\normalsize 3. The Farey Subsequence $\mathcal{F}\bigl(\mathbb{B}(n),m\bigr)$}

Definition~(\ref{eq:2}) implies that a Farey subsequence $\mathcal{F}\bigl(\mathbb{B}(n),m\bigr)$ can be
regarded as the intersection
\begin{equation*}
\mathcal{F}\bigl(\mathbb{B}(n),m\bigr)=\mathcal{F}_n^m\cap\mathcal{G}_n^m\ .
\end{equation*}
Its halfsequences can be described with the help of the following statement:

\begin{theorem} Consider a Farey subsequence $\mathcal{F}\bigl(\mathbb{B}(n),m\bigr)$.
\label{prop:1} The maps
\begin{align}
\label{eq:4} \mathcal{F}^{\leq\frac{1}{2}}\bigl(\mathbb{B}(n),m\bigr)&\to\mathcal{F}^m_{n-m}\ , &
\tfrac{h}{k}&\mapsto \tfrac{h}{k-h}\ , & \left[\begin{smallmatrix}\!h\!\\
\!k\!\end{smallmatrix}\right]&\mapsto\left[\begin{smallmatrix}1&0\\-1&1\end{smallmatrix}\right]\cdot
\left[\begin{smallmatrix}\!h\!\\ \!k\!\end{smallmatrix}\right]\ ,\\ \label{eq:5}
\mathcal{F}^m_{n-m}&\to\mathcal{F}^{\leq\frac{1}{2}}\bigl(\mathbb{B}(n),m\bigr)\ , & \tfrac{h}{k}&\mapsto
\tfrac{h}{k+h}\ , & \left[\begin{smallmatrix}\!h\!\\
\!k\!\end{smallmatrix}\right]&\mapsto\left[\begin{smallmatrix}1&0\\1&1\end{smallmatrix}\right]\cdot
\left[\begin{smallmatrix}\!h\!\\ \!k\!\end{smallmatrix}\right]\ ,\\ \label{eq:12}
\mathcal{F}^{\geq\frac{1}{2}}\bigl(\mathbb{B}(n),m\bigr)&\to\mathcal{G}_m^{2m-n}\ , & \tfrac{h}{k}&\mapsto
\tfrac{2h-k}{h}\ , & \left[\begin{smallmatrix}\!h\!\\
\!k\!\end{smallmatrix}\right]&\mapsto\left[\begin{smallmatrix}2&-1\\1&0\end{smallmatrix}\right]\cdot
\left[\begin{smallmatrix}\!h\!\\ \!k\!\end{smallmatrix}\right]\ , \intertext{and} \label{eq:13}
\mathcal{G}_m^{2m-n}&\to\mathcal{F}^{\geq\frac{1}{2}}\bigl(\mathbb{B}(n),m\bigr)\ , & \tfrac{h}{k}&\mapsto
\tfrac{k}{2k-h}\ , & \left[\begin{smallmatrix}\!h\!\\
\!k\!\end{smallmatrix}\right]&\mapsto\left[\begin{smallmatrix}0&1\\-1&2\end{smallmatrix}\right]\cdot
\left[\begin{smallmatrix}\!h\!\\ \!k\!\end{smallmatrix}\right]\ ,
\end{align}
are order-preserving and bijective.

The maps
\begin{align}
\nonumber \mathcal{F}^{\leq\frac{1}{2}}\bigl(\mathbb{B}(n),m\bigr)&\to\mathcal{G}^{n-2m}_{n-m}\ , &
\tfrac{h}{k}&\mapsto\tfrac{k-2h}{k-h}\ , & \left[\begin{smallmatrix}\!h\!\\
\!k\!\end{smallmatrix}\right]&\mapsto\left[\begin{smallmatrix}-2&1\\-1&1\end{smallmatrix}\right]\cdot
\left[\begin{smallmatrix}\!h\!\\ \!k\!\end{smallmatrix}\right]\ , \\ \nonumber
\mathcal{G}^{n-2m}_{n-m}&\to\mathcal{F}^{\leq\frac{1}{2}}\bigl(\mathbb{B}(n),m\bigr)\ , &
\tfrac{h}{k}&\mapsto\tfrac{k-h}{2k-h}\ , & \left[\begin{smallmatrix}\!h\!\\
\!k\!\end{smallmatrix}\right]&\mapsto\left[\begin{smallmatrix}-1&1\\-1&2\end{smallmatrix}\right]\cdot
\left[\begin{smallmatrix}\!h\!\\ \!k\!\end{smallmatrix}\right]\ , \\ \label{eq:6}
\mathcal{F}^{\geq\frac{1}{2}}\bigl(\mathbb{B}(n),m\bigr)&\to\mathcal{F}^{n-m}_m\ , &
\tfrac{h}{k}&\mapsto\tfrac{k-h}{h}\ , & \left[\begin{smallmatrix}\!h\!\\
\!k\!\end{smallmatrix}\right]&\mapsto\left[\begin{smallmatrix}-1&1\\1&0\end{smallmatrix}\right]\cdot
\left[\begin{smallmatrix}\!h\!\\ \!k\!\end{smallmatrix}\right]\ ,\intertext{and} \label{eq:7}
\mathcal{F}^{n-m}_m&\to\mathcal{F}^{\geq\frac{1}{2}}\bigl(\mathbb{B}(n),m\bigr)\ , &
\tfrac{h}{k}&\mapsto\tfrac{k}{k+h}\ , & \left[\begin{smallmatrix}\!h\!\\
\!k\!\end{smallmatrix}\right]&\mapsto\left[\begin{smallmatrix}0&1\\1&1\end{smallmatrix}\right]\cdot
\left[\begin{smallmatrix}\!h\!\\ \!k\!\end{smallmatrix}\right]\ ,
\end{align}
are order-reversing and bijective.
\end{theorem}

\noindent{\it Proof}. For any integer $h$, $1\leq h\leq m$, we have
\begin{multline*}
\left|\left\{\tfrac{h}{k}\in\mathcal{F}\bigl(\mathbb{B}(n),m\bigr):\
\tfrac{h}{k}<\tfrac{1}{2}\right\}\right|=\phi(h;[2h+1,h+n-m])\\=\sum_{d\in[1,h]:\ d|h}\overline{\mu}(d)\cdot
\left(\left\lfloor\tfrac{h+n-m}{d}\right\rfloor-\tfrac{2h}{d}\right)=\sum_{d\in[1,h]:\
d|h}\overline{\mu}(d)\cdot
\left(\left\lfloor\tfrac{n-m}{d}\right\rfloor-\tfrac{h}{d}\right)\\=\phi(h;[h+1,n-m])=\left|\left\{\tfrac{h}{k}\in
\mathcal{F}_{n-m}^m:\ \tfrac{h}{k}<\tfrac{1}{1}\right\}\right|
\end{multline*}
and see that the sequences $\mathcal{F}^{\leq\frac{1}{2}}\bigl(\mathbb{B}(n),m\bigr)$ and
$\mathcal{F}_{n-m}^m$ are of the same cardinality, but this conclusion also implies
$\bigl|\mathcal{F}^{\geq\frac{1}{2}}\bigl(\mathbb{B}(n),m\bigr)\bigr|=|\mathcal{F}^{n-m}_m|$, due to
bijection~(\ref{eq:3}). The proof of the assertions concerning maps~(\ref{eq:4}), (\ref{eq:5}), (\ref{eq:6})
and~(\ref{eq:7}) is completed by checking that, on the one hand, a fraction $\tfrac{h_j}{k_j}$ precedes a
fraction $\tfrac{h_{j+1}}{k_{j+1}}$ in $\mathcal{F}_{n-m}^m$ if and only if the fraction
$\tfrac{h_j}{k_j+h_j}$ precedes the fraction $\tfrac{h_{j+1}}{k_{j+1}+h_{j+1}}$ in
$\mathcal{F}^{\leq\frac{1}{2}}\bigl(\mathbb{B}(n),m\bigr)$; on the other hand, a fraction $\tfrac{h_j}{k_j}$
precedes a fraction $\tfrac{h_{j+1}}{k_{j+1}}$ in $\mathcal{F}^{n-m}_m$ if and only if
$\tfrac{k_{j+1}}{k_{j+1}+h_{j+1}}$ precedes $\tfrac{k_j}{k_j+h_j}$
in~$\mathcal{F}^{\geq\frac{1}{2}}\bigl(\mathbb{B}(n),m\bigr)$. The remaining assertions of the theorem now
follow, thanks to Lemma~\ref{prop:2}. \hfill$\Box$

For example,
\begin{align*}
\mathcal{F}_6 &= \bigl(\tfrac{0}{1}<\tfrac{1}{6}<\tfrac{1}{5}<\tfrac{1}{4}<\tfrac{1}{3}<\tfrac{2}{5}<
\tfrac{1}{2}<\tfrac{3}{5}<\tfrac{2}{3}<\tfrac{3}{4}<\tfrac{4}{5}<\tfrac{5}{6}<\tfrac{1}{1}\bigr)\ ,\\
\mathcal{F}_6^4 &= \bigl(\tfrac{0}{1}<\tfrac{1}{6}<\tfrac{1}{5}<\tfrac{1}{4}<\tfrac{1}{3}<\tfrac{2}{5}<
\tfrac{1}{2}<\tfrac{3}{5}<\tfrac{2}{3}<\tfrac{3}{4}<\tfrac{4}{5} \phantom{<}\,\,\, \phantom{\tfrac{5}{6}}<
\tfrac{1}{1}\bigr)\ ,\\ \mathcal{G}_6^4 &=
\bigl(\tfrac{0}{1}\phantom{<}\,\,\,\phantom{\tfrac{1}{6}}\phantom{<}\,\,\,\phantom{\tfrac{1}{5}}
\phantom{<}\,\,\,\phantom{\tfrac{1}{4}}\,<
\tfrac{1}{3}\phantom{<}\,\,\,\phantom{\tfrac{2}{5}}<\tfrac{1}{2}<\tfrac{3}{5}<\tfrac{2}{3}<\tfrac{3}{4}<
\tfrac{4}{5}< \tfrac{5}{6}<\tfrac{1}{1}\bigr)\ ,\\
\mathcal{F}\bigl(\mathbb{B}(6),4\bigr)&=\bigl(\tfrac{0}{1}\phantom{<}\,\,\,
\phantom{\tfrac{1}{6}}\phantom{<}\,\,\,\phantom{\tfrac{1}{5}}\phantom{<}\,\,\,
\phantom{\tfrac{1}{4}}\,<\tfrac{1}{3}\phantom{<}\,\,\,\phantom{\tfrac{2}{5}}<
\tfrac{1}{2}<\tfrac{3}{5}<\tfrac{2}{3}<\tfrac{3}{4}<\tfrac{4}{5}\phantom{<}\,\,\,
\phantom{\tfrac{5}{6}}<\tfrac{1}{1}\bigr)\ ,\\
\mathcal{F}_{6-4}^4=\mathcal{F}_2&=\bigl(\tfrac{0}{1}\phantom{<}\,\,\,
\phantom{\tfrac{1}{6}}\phantom{<}\,\,\,\phantom{\tfrac{1}{5}}\phantom{<}\,\,\,
\phantom{\tfrac{1}{4}}\,<\tfrac{1}{2}\phantom{<}\,\,\,\phantom{\tfrac{2}{5}} < \tfrac{1}{1}\bigr)\ ,\\
\mathcal{G}_4^{2\cdot 4-6}=\mathcal{G}_4^2 &=\phantom{\bigl(}\phantom{\tfrac{0}{1}}\phantom{<}\,\,\,
\phantom{\tfrac{1}{6}}\phantom{<}\,\,\,\phantom{\tfrac{1}{5}}\phantom{<}\,\,\,
\phantom{\tfrac{1}{4}}\,\phantom{<}\phantom{\tfrac{1}{3}}\phantom{<}\,\,\,\phantom{\tfrac{2}{5}}\,\,\,\,
\,\,\,\,\,\bigl(\tfrac{0}{1}<\tfrac{1}{3}<\tfrac{1}{2}<\tfrac{2}{3}<\tfrac{3}{4}\phantom{<}
\phantom{\tfrac{5}{6}}\,\,\,<\tfrac{1}{1}\bigr)\ ,\\ \mathcal{G}_{6-4}^{6-2\cdot 4}
=\mathcal{F}_2&=\bigl(\tfrac{0}{1}\phantom{<}\,\,\,
\phantom{\tfrac{1}{6}}\phantom{<}\,\,\,\phantom{\tfrac{1}{5}}\phantom{<}\,\,\,
\phantom{\tfrac{1}{4}}\,<\tfrac{1}{2}\phantom{<}\,\,\,\phantom{\tfrac{2}{5}} < \tfrac{1}{1}\bigr)\ ,\\
\mathcal{F}_4^{6-4}=\mathcal{F}_4^2&=\phantom{\bigl(}\phantom{\tfrac{0}{1}}\phantom{<}\,\,\,
\phantom{\tfrac{1}{6}}\phantom{<}\,\,\,\phantom{\tfrac{1}{5}}\phantom{<}\,\,\,
\phantom{\tfrac{1}{4}}\,\phantom{<}\phantom{\tfrac{1}{3}}\phantom{<}\,\,\,\phantom{\tfrac{2}{5}}\,\,\,\,
\,\,\,\,\,\bigl(\tfrac{0}{1}<\tfrac{1}{4}<\tfrac{1}{3}<\tfrac{1}{2}<\tfrac{2}{3}\phantom{<}
\phantom{\tfrac{5}{6}}\,\,\,<\tfrac{1}{1}\bigr)\ .
\end{align*}

Theorem~\ref{prop:1} allows us to write down the formula
\begin{equation*}
\left|\mathcal{F}\bigl(\mathbb{B}(n),m\bigr)\right|=|\mathcal{F}_{n-m}^m|+|\mathcal{F}^{n-m}_m|-1=
\bigl|\mathcal{F}_{n-m}^{\min\{m,n-m\}}\bigr|+\bigl|\mathcal{F}^{\min\{m,n-m\}}_m\bigr|-1\ ,
\end{equation*}
for the number of elements of a sequence $\mathcal{F}\bigl(\mathbb{B}(n),m\bigr)$, cf.~[10,
Proposition~7.3(ii)]. Recall that $|\mathcal{F}_q^p|=1+\sum_{d\geq
1}\overline{\mu}(d)\cdot\left(\left\lfloor\tfrac{q}{d}\right\rfloor-
\tfrac{1}{2}\left\lfloor\tfrac{p}{d}\right\rfloor\right)\cdot\left\lfloor\tfrac{p}{d}+1\right\rfloor$, for
any Farey subsequence $\mathcal{F}_q^p$ with $0<p\leq q$, see~[10, Remark~7.10(ii)(b)]; notice that the
well-known relation $\sum_{d\geq 1}\overline{\mu}(d)\cdot\left\lfloor\tfrac{t}{d}\right\rfloor=1$, for any
positive integer $t$, leads us to one more formula: $|\mathcal{F}_q^p|=\tfrac{3}{2}+\sum_{d\geq
1}\overline{\mu}(d)\cdot\left\lfloor\tfrac{p}{d}\right\rfloor\cdot\left(\left\lfloor\tfrac{q}{d}\right\rfloor-
\tfrac{1}{2}\left\lfloor\tfrac{p}{d}\right\rfloor\right)$.

Thus, we have
\begin{multline*}
\left|\mathcal{F}\bigl(\mathbb{B}(n),m\bigr)\right|=\frac{3}{2}+\sum_{d\geq 1}\overline{\mu}(d)\!\cdot\!
\left\lfloor\tfrac{\min\{m,n-m\}}{d}\right\rfloor\!\cdot\!
\left(\!\left\lfloor\tfrac{n-m}{d}\right\rfloor\!-\tfrac{1}{2}\!\left\lfloor\tfrac{\min\{m,n-m\}}{d}
\right\rfloor\right)
\\
+\frac{3}{2}+\sum_{d\geq 1}\overline{\mu}(d)\!\cdot\!
\left\lfloor\tfrac{\min\{m,n-m\}}{d}\right\rfloor\!\cdot\!
\left(\!\left\lfloor\tfrac{m}{d}\right\rfloor\!-\tfrac{1}{2}\!\left\lfloor\tfrac{\min\{m,n-m\}}{d}
\right\rfloor\right)-1
\end{multline*}
or
\begin{equation*}
\left|\mathcal{F}\bigl(\mathbb{B}(n),m\bigr)\right|=2+\sum_{d\geq 1}\overline{\mu}(d)\!\cdot\!
\left\lfloor\tfrac{\min\{m,n-m\}}{d}\right\rfloor\!\cdot\!
\left(\!\left\lfloor\tfrac{n-m}{d}\right\rfloor\!+\!\left\lfloor\tfrac{m}{d}\right\rfloor\!-\!
\left\lfloor\tfrac{\min\{m,n-m\}}{d}\right\rfloor\!\right)\ ,
\end{equation*}
that is,
\begin{equation*}
\left|\mathcal{F}\bigl(\mathbb{B}(n),m\bigr)\right|-2=\sum_{d\geq
1}\overline{\mu}(d)\!\cdot\!\left\lfloor\tfrac{m}{d}\right\rfloor\!\!\cdot\!\!\left\lfloor\tfrac{n-m}{d}
\right\rfloor\ .
\end{equation*}
In particular, for any positive integer $t$ we have
\begin{equation*}
\sum_{d\geq 1}\overline{\mu}(d)\!\cdot\!\left\lfloor\tfrac{t}{d}\right\rfloor^2=
\left|\mathcal{F}\bigl(\mathbb{B}(2t),t\bigr)\right|-2=2\left|\mathcal{F}_t\right|-3\ .
\end{equation*}

\begin{proposition}
Consider a Farey subsequence $\mathcal{F}\bigl(\mathbb{B}(n),m\bigr)$.
\begin{itemize}
\item[\rm(i)] Suppose $m\geq\frac{n}{2}$. The map
\begin{align}
\label{eq:9} \mathcal{F}^{\leq\frac{1}{2}}\bigl(\mathbb{B}(n),m\bigr)&\to
\mathcal{F}^{\leq\frac{1}{2}}\bigl(\mathbb{B}(n),m\bigr)\ , & \tfrac{h}{k}&\mapsto\tfrac{k-2h}{2k-3h}\ , &
\left[\begin{smallmatrix}\!h\!\\
\!k\!\end{smallmatrix}\right]&\mapsto\left[\begin{smallmatrix}-2&1\\-3&2\end{smallmatrix}\right]\cdot
\left[\begin{smallmatrix}\!h\!\\ \!k\!\end{smallmatrix}\right]\ , \intertext{is order-reversing and
bijective. The map}\nonumber \mathcal{F}^{\leq\frac{1}{2}}\bigl(\mathbb{B}(n),m\bigr)&\to
\mathcal{F}^{\geq\frac{1}{2}}\bigl(\mathbb{B}(n),m\bigr)\ , & \tfrac{h}{k}&\mapsto\tfrac{k-h}{2k-3h}\ , &
\left[\begin{smallmatrix}\!h\!\\
\!k\!\end{smallmatrix}\right]&\mapsto\left[\begin{smallmatrix}-1&1\\-3&2\end{smallmatrix}\right]\cdot
\left[\begin{smallmatrix}\!h\!\\ \!k\!\end{smallmatrix}\right]\ , \intertext{is order-preserving and
injective. The map}\nonumber \mathcal{F}^{\leq\frac{1}{2}}\bigl(\mathbb{B}(n),m\bigr)&\to
\mathcal{F}^{\geq\frac{1}{2}}\bigl(\mathbb{B}(n),m\bigr)\ , & \tfrac{h}{k}&\mapsto\tfrac{k-h}{k}\ , &
\left[\begin{smallmatrix}\!h\!\\
\!k\!\end{smallmatrix}\right]&\mapsto\left[\begin{smallmatrix}-1&1\\0&1\end{smallmatrix}\right]\cdot
\left[\begin{smallmatrix}\!h\!\\ \!k\!\end{smallmatrix}\right]\ ,
\end{align}
is order-reversing and injective.
\item[\rm(ii)]
Suppose $m\leq\frac{n}{2}$. The map
\begin{align}
\label{eq:10} \mathcal{F}^{\geq\frac{1}{2}}\bigl(\mathbb{B}(n),m\bigr)&\to
\mathcal{F}^{\geq\frac{1}{2}}\bigl(\mathbb{B}(n),m\bigr)\ , & \tfrac{h}{k}&\mapsto\tfrac{h}{3h-k}\ , &
\left[\begin{smallmatrix}\!h\!\\
\!k\!\end{smallmatrix}\right]&\mapsto\left[\begin{smallmatrix}1&0\\3&-1\end{smallmatrix}\right]\cdot
\left[\begin{smallmatrix}\!h\!\\ \!k\!\end{smallmatrix}\right]\ , \intertext{is order-reversing and
bijective. The map}\nonumber \mathcal{F}^{\geq\frac{1}{2}}\bigl(\mathbb{B}(n),m\bigr)&\to
\mathcal{F}^{\leq\frac{1}{2}}\bigl(\mathbb{B}(n),m\bigr)\ , & \tfrac{h}{k}&\mapsto\tfrac{2h-k}{3h-k}\ , &
\left[\begin{smallmatrix}\!h\!\\
\!k\!\end{smallmatrix}\right]&\mapsto\left[\begin{smallmatrix}2&-1\\3&-1\end{smallmatrix}\right]\cdot
\left[\begin{smallmatrix}\!h\!\\ \!k\!\end{smallmatrix}\right]\ , \intertext{is order-preserving and
injective. The map}\nonumber \mathcal{F}^{\geq\frac{1}{2}}\bigl(\mathbb{B}(n),m\bigr)&\to
\mathcal{F}^{\leq\frac{1}{2}}\bigl(\mathbb{B}(n),m\bigr)\ , & \tfrac{h}{k}&\mapsto\tfrac{k-h}{k}\ , &
\left[\begin{smallmatrix}\!h\!\\
\!k\!\end{smallmatrix}\right]&\mapsto\left[\begin{smallmatrix}-1&1\\0&1\end{smallmatrix}\right]\cdot
\left[\begin{smallmatrix}\!h\!\\ \!k\!\end{smallmatrix}\right]\ ,
\end{align}
is order-reversing and injective.
\end{itemize}
\end{proposition}

\noindent{\it Proof}. To prove that map~(\ref{eq:9}) is order-reversing and bijective, notice that
$\mathcal{F}_{n-m}^m= \mathcal{F}_{n-m}$, and consider the composite map
\begin{equation*}
\tfrac{h}{k}\ \stackrel{
{\mathcal{F}^{\leq\frac{1}{2}}(\mathbb{B}(n),m)\xrightarrow{\rm(\ref{eq:4})}\mathcal{F}_{n-m}}}{\mapsto}\
\tfrac{h}{k-h}\ \stackrel{ {\mathcal{F}_{n-m}\xrightarrow{\rm(\ref{eq:11})}\mathcal{F}_{n-m}}}{\mapsto}\
\tfrac{k-2h}{k-h}\ \stackrel{\mathcal{F}_{n-m}
{\xrightarrow{\rm(\ref{eq:5})}\mathcal{F}^{\leq\frac{1}{2}}(\mathbb{B}(n),m)}}{\mapsto}\ \tfrac{k-2h}{2k-3h}\
.
\end{equation*}
Similarly, under the hypothesis of assertion~(ii) we have $\mathcal{G}_{m}^{2m-n}=\mathcal{F}_{m}$; consider
the composite map
\begin{equation*}
\tfrac{h}{k}\ \stackrel{
{\mathcal{F}^{\geq\frac{1}{2}}(\mathbb{B}(n),m)\xrightarrow{\rm(\ref{eq:12})}\mathcal{F}_m}}{\mapsto}\
\tfrac{2h-k}{h}\ \stackrel{ {\mathcal{F}_m\xrightarrow{\rm(\ref{eq:11})}\mathcal{F}_m}}{\mapsto}\
\tfrac{k-h}{h}\ \stackrel{\mathcal{F}_m
{\xrightarrow{\rm(\ref{eq:13})}\mathcal{F}^{\geq\frac{1}{2}}(\mathbb{B}(n),m)}}{\mapsto}\ \tfrac{h}{3h-k}
\end{equation*}
to see that map~(\ref{eq:10}) is order-reversing and bijective.

The remaining assertions follow from the observation that map~(\ref{eq:3}) is the order-reversing
bijection.\hfill $\Box$

We conclude the note by listing a few pairs of fractions adjacent in
$\mathcal{F}\bigl(\mathbb{B}(n),m\bigr)$, \mbox{see~[8, 9]} on such pairs within the sequences
$\mathcal{F}\bigl(\mathbb{B}(2m),m\bigr)$; we find the neighbors of the images of several fractions of
$\mathcal{F}\bigl(\mathbb{B}(n),m\bigr)$ under bijections from Theorem~\ref{prop:1}, and we then reflect them
back to $\mathcal{F}\bigl(\mathbb{B}(n),m\bigr)$:

\begin{remark}
Consider a Farey subsequence $\mathcal{F}\bigl(\mathbb{B}(n),m\bigr)$, with $n\neq 2m$.
\begin{itemize}
\item[\rm(i)]
Suppose $m>\frac{n}{2}$.

The fraction $\tfrac{n-m-1}{2(n-m)-1}$ precedes~$\tfrac{1}{2}$, and the fraction $\tfrac{n-m+1}{2(n-m)+1}$
succeeds~$\tfrac{1}{2}$.

The fraction $\tfrac{2\min\{n-m,\lfloor\frac{m+1}{2}\rfloor\}-1}{3\min\{n-m,\lfloor\frac{m+1}{2}\rfloor\}-1}$
precedes $\tfrac{2}{3}$; the fraction
$\tfrac{2\min\{n-m,\lfloor\frac{m-1}{2}\rfloor\}+1}{3\min\{n-m,\lfloor\frac{m-1}{2}\rfloor\}+1}$ succeeds
$\tfrac{2}{3}$.

If we additionally have $n-m>1$, then the fraction
\begin{equation*}
\begin{cases}\frac{n-m-2}{2}\!\Bigm/\!\frac{3(n-m)-4}{2}, &\text{\em if $n-m$ is even},\\
\frac{n-m-1}{2}\!\Bigm/\!\frac{3(n-m)-1}{2}, &\text{\em if $n-m$ is odd},
\end{cases}
\end{equation*}
precedes $\tfrac{1}{3}$; the fraction
\begin{equation*}
\begin{cases}\frac{n-m}{2}\!\Bigm/\!\frac{3(n-m)-2}{2}, &\text{\em if $n-m$ is even},\\
\frac{n-m+1}{2}\!\Bigm/\!\frac{3(n-m)+1}{2}, &\text{\em if $n-m$ is odd},
\end{cases}
\end{equation*}
succeeds $\tfrac{1}{3}$.

\item[\rm(ii)]
Suppose $m<\frac{n}{2}$.

The fraction $\tfrac{m}{2m+1}$ precedes~$\tfrac{1}{2}$, and the fraction $\tfrac{m}{2m-1}$ succeeds
$\tfrac{1}{2}$.

The fraction $\tfrac{\min\{m,\lfloor\frac{n-m+1}{2}\rfloor\}-1}{3\min\{m,\lfloor\frac{n-m+1}{2}\rfloor\}-2}$
precedes $\tfrac{1}{3}$; the fraction
$\tfrac{\min\{m,\lfloor\frac{n-m-1}{2}\rfloor\}+1}{3\min\{m,\lfloor\frac{n-m-1}{2}\rfloor\}+2}$ succeeds
$\tfrac{1}{3}$.

If we additionally have $m>1$, then the fraction
\begin{equation*}
\begin{cases}(m-1)\!\Bigm/\!\frac{3m-2}{2}, &\text{\em if $m$ is even},\\
m\!\Bigm/\!\frac{3m+1}{2}, &\text{\em if $m$ is odd},
\end{cases}
\end{equation*}
precedes $\tfrac{2}{3}$; the fraction
\begin{equation*}
\begin{cases}(m-1)\!\Bigm/\!\frac{3m-4}{2}, &\text{\em if $m$ is even},\\
m\!\Bigm/\!\frac{3m-1}{2}, &\text{\em if $m$ is odd},
\end{cases}
\end{equation*}
succeeds $\tfrac{2}{3}$.
\end{itemize}
\end{remark}

\vskip 30pt

\section*{\normalsize References}

\noindent[1] D.~Acketa and J.~\v{Z}uni\'{c}, {\em On the number of linear partitions of the $(m,n)$-grid},
Inform. Process. Lett. {\bf 38} (1991), no.~3, 163--168.

\noindent[2] M.~Aigner, {\em Combinatorial Theory, Reprint of the 1979 original}, Classics in Mathematics,
Springer-Verlag, Berlin, 1997.

\noindent[3] A.A.~Buchstab, {\em Teoriya Chisel} (in Russian) [{\em Number Theory}], Uchpedgiz, Moscow, 1960.

\noindent[4] S.B.~Gashkov and V.N.~Chubarikov, {\em Arifmetika. Algoritmy. Slozhnost' Vychisleniy, Third
edition} (in Russian) [{\em Arithmetic. Algorithms. Complexity of Computation}], Drofa, Moscow, 2005.

\noindent[5] R.L.~Graham, D.E.~Knuth and O.~Patashnik, {\em Concrete Mathematics. A Foundation for Computer
Science, Second edition}, Addison-Wesley, Reading Massachusetts, 1994.

\noindent[6] G.H.~Hardy and E.M.~Wright, {\em An Introduction to the Theory of Numbers, Fifth edition},
Clarendon Press, Oxford, 1979.

\noindent[7] A.O.~Matveev, {\em A note on Boolean lattices and Farey sequences}, Integers {\bf 7} (2007),
A20.

\noindent[8] A.O.~Matveev, {\em Neighboring fractions in Farey subsequences}, {\tt arXiv:0801.1981} .

\noindent[9] A.O.~Matveev, {\em Pattern recognition on oriented matroids: Layers of tope committees}, {\tt
arXiv:math/0612369} .

\noindent[10] A.O.~Matveev, {\em Relative blocking in posets},  J.~Comb.~Optim. {\bf 13} (2007), no.~4,
\mbox{379--403}.

\noindent[11] I.~Niven, H.S.~Zuckerman and H.L.~Montgomery, {\em An Introduction to the Theory of~Numbers,
Fifth edition}, John Wiley \& Sons, Inc., New~York, 1991.

\noindent[12] V.V.~Prasolov, {\em Zadachi po Algebre, Arifmetike i Analizu} (in Russian) [{\em Problems in
Algebra, Arithmetic and Analysis}], Moskovski\u\i\ Tsentr Nepreryvnogo Matematicheskogo Obrazovaniya (MCCME),
Moscow, 2005.

\noindent[13] N.J.A.~Sloane, {\em The On-Line Encyclopedia of Integer Sequences}, published electronically at
{\tt www.research.att.com/\verb"~"\!njas/sequences/}\ , 2008.

\noindent[14] J.J.~Tattersall, {\em Elementary Number Theory in Nine Chapters, Second edition}, Cambridge
University Press, Cambridge, 2005. 
\end{document}